\documentclass[preprint,12pt]{article}

\usepackage[ruled]{algorithm2e}

\usepackage[english]{babel}
\usepackage[utf8x]{inputenc}
\usepackage[T1]{fontenc}
\usepackage{amssymb}
\usepackage{enumitem}
\usepackage{amsthm}
\usepackage{amsmath}
\usepackage[]{algorithm2e}
\usepackage{comment}
\usepackage{thmtools}
\usepackage{lscape}
\usepackage{graphicx}
\usepackage{tabularx}
\usepackage{amsmath, amssymb}
\usepackage[T1]{fontenc}
\usepackage{url}
\usepackage{subfig}
\usepackage{xcolor}

\usepackage{alltt}
\usepackage{adjustbox}
\usepackage{cite}
\usepackage{verbatim}
\usepackage{xcolor}
\usepackage[all]{xy}
\usepackage{tikz}
\usepackage{multirow}
\usepackage{siunitx}
\usepackage{dcolumn}

\usepackage{hyperref}

\newcolumntype{d}[1]{D{.}{.}{#1}}

\renewcommand{\geq}{\geqslant}
\renewcommand{\leq}{\leqslant}

\usepackage{epsfig}
\usepackage{makeidx}         % allows index generation
\usepackage{graphicx}        % standard LaTeX graphics tool
\usepackage{url}
\renewcommand{\Re}{{\mathbb R}}

\newcommand{\vect} \mathbf

\newcommand{\glebdone}[1]{}

\newcommand{\wudone}[1]{}

\newcommand{\peiqidone}[1]{}

\makeatletter
\newcommand*{\itemequation}[3][]{%
  \item
  \begingroup
    \refstepcounter{equation}%
    \ifx\\#1\\%
    \else  
      \label{#1}%
    \fi
    \sbox0{#2}%
    \sbox2{$\displaystyle#3\m@th$}%
    \sbox4{\@eqnnum}%
    \dimen@=.5\dimexpr\linewidth-\wd2\relax
    \ifcase
        \ifdim\wd0>\dimen@
          \z@
        \else
          \ifdim\wd4>\dimen@
            \z@
          \else 
            \@ne
          \fi 
        \fi
      \@latex@warning{Equation is too large}%
    \fi
    \noindent   
    \rlap{\copy0}%
    \rlap{\hbox to \linewidth{\hfill\copy2\hfill}}%
    \hbox to \linewidth{\hfill\copy4}%
    \hspace{0pt}% allow linebreak
  \endgroup
  \ignorespaces 
}
\makeatother

\usepackage{authblk}

\title{Fitting scattered data with optional monotonicity constraints on GPU: LipFit package}

\author[1]{Gleb Beliakov}

\affil[1]{School of Information Technology, Deakin University, Burwood 3125, Australia}
\begin{document}

\maketitle

\abstract{
This paper presents a method of multivariate scattered data interpolation and approximation that produces optimal Lipschitz-continuous approximation, subject to the desired monotonicity constraints. This method relies on tight upper and lower approximations to the data, and is similar in its spirit to the nearest-neighbour approximation but does not suffer from discontinuities. Local Lipschitz interpolation and Lipschitz smoothing are also presented. This approach falls under the umbrella of instance-based approximation with no training phase, and it is suitable for GPU-based parallelisation. A Python GPU-friendly package LipFit which implements the methods discussed  is discussed.
}


\section{Introduction}

Throughout this text $d$ will denote the dimensionality of the
space, and $N$ will denote the size of the data set. We are given a
data set representing the values of an unknown function $f$ in tabular format.

\begin{center}
\begin{tabular}{|r|r|r|r|r|}\hline
$x_1$  & $x_2$ & $x_3$ &  $x_4$ & $y$ \\
 \hline
$x_1^1$  & $x_2^1$ & $x_3^1$ &  $x_4^1$ & $y^1$ \\
$x_1^2$  & $x_2^2$ & $x_3^2$ &  $x_4^2$ & $y^2$ \\
$x_1^3$  & $x_2^3$ & $x_3^3$ &  $x_4^3$ & $y^3$ \\
$\vdots$ & & & & \\
$x_1^N$  & $x_2^N$ & $x_3^N$ &  $x_4^N$ & $y^N$ \\
\hline
\end{tabular}
\end{center}

There is no special structure in the data set, i.e., the data are \textit{scattered}.
We assume that the data set $D$ was generated by a function $f$ (which we call the underlying function, it is unknown to us) which satisfies Lipschitz condition with
a Lipschitz constant $M$:
$$
 |f(x)-f(z)|\leq M d(x,z),
$$
for all $x$ and $z$, where $d(x,z)$ is a distance function.
We look for a fitting function  $g \approx f$, such that
$$
g(x^k)=y^k, k=1,\ldots,N, \mbox{ (interpolation conditions)}
$$
which provides the best uniform approximation to $f$ in the worst case scenario, i.e.,
$g$ minimizes the maximal possible error at any $x$
$$
\max_{f} \max_{x \in X} |f(x)-g(x)|.
$$
The maximum is taken over all possible Lipschitz functions that could have produced the data set $D$.

When the interpolation conditions are inconsistent with the Lipschitz properties (no Lipschitz function with constant $M$ could have produced $D$), or when interpolation is not required (e.g., noisy data), we aim at approximating the data
$$
g(x^k)\approx y^k, k=1,\ldots,N, \mbox{ (approximation conditions)}
$$
which can be understood in the least squares, least absolute deviations or another appropriate sense.

In addition, we enforce various types of monotonicity constraints on $g$, coordinatewise, which narrows down the range of $f$ and hence $g$.

\section{Literature review}

\subsection{Earlier works}

Multivariate data interpolation and approximation is a very common
problem in many branches of science. There is a great number of
techniques developed for various instances of this problem, such as
polynomial regression, spline interpolation and smoothing, wavelets,
nearest neighbour search, Sibson interpolation, MARS (multivariate
adaptive regression splines), machine learning techniques (e.g.,
decision trees), neural networks, radial basis functions, etc. For
an overview the reader is referred to
\cite{Alfeld1989_inbook,Cherkassky1998_book, Hastie2001_book}.

Monotonicity with respect to one or more variables is one such
property which frequently arises in practice. There are many tasks
that require monotone interpolation and approximation in several
variables. Here we mention just five examples coming from different
sciences.
\begin{itemize}
    \item Dose-response curves and surfaces in biochemistry and pharmacology;
    \item Design of aggregation operators in multicriteria decision making and fuzzy logic
    \cite{Beliakov2007_book};
    \item Approximation of copulas and quasi-copulas in statistics;
    \item Empirical option pricing models \cite{Hutchison1994_JF} in finance;
    \item Approximation of potential functions in physical and chemical systems.
\end{itemize}

In the
bivariate case, monotonicity preserving splines were studied in
\cite{Beatson1985_SIAM,Carlson1985_SIAM,Costantini1990_SIAM,Costantini1996_CAGD1}.
In these methods the data should be typically given on a rectangular
grid. The case of bivariate scattered data with tensor product
topology is discussed in
\cite{Costantini1996_CAGD,Costantini1996_JCAM1}. Only few methods
are able to deal with scattered data
\cite{Utreras1991_CA,Han1997_SIAM,Costantini1999_CAGD,Costantini1996_JCAM},
which is typical in applications.

Monotone tensor product regression splines were used for more than
two variables in \cite{Beliakov2000_ata,Beliakov2002_ijfks},
including the case of scattered data. However tensor product
schemata have significant drawbacks in the multivariate setting,
which stem from the fact that an exponential number of basis
functions (and hence spline coefficients to compute) is required.
With the increasing dimension, this number quickly exceeds the
number of data, leading to ill-conditioned systems of equations. For
more than five variables tensor product splines are not practical.

Triangulation based monotone splines are suitable for dealing with
scattered data. In the bivariate case they were studied in
\cite{Willemans1996_NA,Costantini1999_CAGD,Costantini1996_JCAM}. One
drawback of triangulation based schemata is the lack of continuous
dependency of the interpolant on the data \cite{Alfeld1989_inbook}.
Small changes in the abscissae of data points may lead to a
completely different triangulation, which will drastically change
the behavior of the spline. Furthermore, building triangulation in
more than two variables is computationally expensive; it becomes
prohibitive for more than 10 variables. The number of elements of
such triangulation also grows exponentially with the dimension.
Splines on triangulations have also another drawback. Even if the
data is monotone, it is not always possible to build a monotone
interpolating spline. 

There are alternative methods of multivariate interpolation and
approximation of scattered data, like radial basis functions,
$k$-nearest neighbors approximation, Sibson's natural neighbor
interpolation, and neural networks. These methods do not incorporate
monotonicity, although some attempts to enforce monotonicity of
kernel based approximation (in the univariate case) and of the
neural networks have been made \cite{Hall2001_AnnStat,
Hutchison1994_JF}.

\subsection{Recent works}

Machine learning approaches include neural networks, decision trees, $k$-nearest neighbour approximation and various hybrid methods. Twin  neural networks combined with kNN methods are presented in \cite{wetzel2023twinneuralnetworkimproved}. Random kernel kNN use ensemble techniques to improve on traditional kNN \cite{10.3389/fdata.2024.1402384}. Another recent work involving kNN  \cite{koutensky2024overcoming,koutensky2022fast}  focuses on improving inference time. The work \cite{kravberg2022activenearestneighborregression} looks at natural neighbours approaches through Delaunay refinements.

A recent survey of shape-preserving approximation in one dimension is presented in \cite{SMAI-JCM_2019__S5__99_0}. Reference \cite{ doi:10.1142/S0219530521500299} discusses  characterisation and continuity of multivariate monotone regression.
The monotonicity preserving approximation has a number of challenges related to the curse of dimensionality, reviewed in \cite{KUNSCH201933}. The interplay between Lipschitz and monotonic approximation is presented in \cite{takatsu2026isotoniclipschitzregressionnew}.

The recent works addressing monotonicity preservation include certified Kolmogorov-Arnold networks \cite{polomolina2024monokan}, monotonic networks \cite{ZHAO2024457}, smooth min-max networks \cite{igel2024smoothminmaxmonotonicnetworks}. Testing monotonicity in ML was discussed in \cite{sharma2020testingmonotonicitymachinelearning} and \cite{liu2022certifiedmonotonicneuralnetworks}.

\section{Lipschitz interpolation}

Lipschitz condition is easy to interpret in terms of the underlying problem. It is simply the upper bound on the rate of change of function
$f$. No differentiability of $f$ is required.

Our goal is to find an interpolant $g$ which approximates $f$ well
at  the points $x$ distinct from the data, given that $f$ is
Lipschitz. We are interested in reliable approximation of $f$, which
means that we want to obtain a good approximation regardless of how
inconvenient $f$ is, even in the worst case scenario.
That is, we solve the following problem.\\

\noindent Find the best interpolating function $g: \Re^d \rightarrow
\Re$,
\begin{equation}
g=\arg \inf\{ \max_{f \in Lip(M) } ||f-g||_{C(X)}, \}
\label{interp_prob}
\end{equation}
such that
$$
g(x^k)=f(x^k)=y^k, k=1,\ldots,N.
$$
$Lip(M)$ denotes the class of functions whose Lipschitz constant is
smaller or equal to $M$. The norm is the uniform norm in the space of continuous functions $C(X)$.

The  method used in \texttt{LipFit} relies on building tight upper
and lower approximations to $f$, denoted by $H^{upper}$ and
$H^{lower}$, proposed in \cite{Beliakov2006_BIT,Beliakov2004_JCAM,Beliakov2007-zd}

\begin{eqnarray}
H^{upper}(x)=\min_k (y^k +M dist(x,x^k)),\nonumber \\
H^{lower}(x) = \max_k (y^k- M dist(x,x^k)).
\label{upperlower}
\end{eqnarray}

Here $dist$ is an appropriately chosen distance in $\Re^d$, in particular the Euclidean distance (other distances that facilitate treatment of monotonicity are presented later on).

Let
\begin{equation}\label{eq:g}
g(x)=\frac{1}{2}(H^{lower}(x)+H^{upper}(x)), \forall x \in X.\end{equation}

Then $g$ is the solution to the Interpolation Problem (\ref{interp_prob}) over the set of all continuous functions
$X \rightarrow R$ that interpolate the data, i.e.,
$$
g = \arg \min_h \max_{f \in Lip(M)} ||f-h||_{C(X)},
$$
and the best error bound is
$$
\max_{f \in Lip(M)} ||f-g||_{C(X)} = M \max_{x \in X}
\min_{k=1,\ldots,N} dist(x,x^k).
$$

Such an interpolation scheme is called \emph{optimal interpolation} and was treated in \cite{Traub1980_book}.

\section{Construction of the interpolant and its features}

Equations (\ref{upperlower}) provide the way to evaluate the
interpolant $g(x)$ at any query point $x$. There is no need for any
preprocessing, and the number of basic arithmetic operations is
proportional to the size of the data set $N$.

It is convenient to view this process from the nearest neighbour approximation perspective (kNN-type methods), although in contrast to kNN, formulas \eqref{upperlower}, \eqref{eq:g} naturally produce a Lipschitz continuous approximation. Consider first the upper bound $H^{upper}$, which is a piecewise linear Lipschitz continuous function with Lipschitz constant $M$. Its evaluation involves choosing the ``nearest'' point $x^k$ which gives the smallest upper bound on the value of any Lipschitz function interpolating $y_k$ at that point. The point $x^k$ is not the nearest in terms of the distance $dist$, but its choice is also influenced by $y_k$. Yet it happens to be among the closest points to $x$.

The lower bound is evaluated similarly, but the choice of the ``nearest'' point $x^k$ may be different to the upper bound. Therefore there are two neighbours of $x$ that influence the value of $g$, and they are usually not on the same side of $x$ (they are on the opposite sides in the one-dimensional case). Hence this method is closer to the natural neighbour approximation \cite{Sibson1981_inbook}, where the neighbours are distributed all around $x$. 

The lower and upper bounds give the maximal range of values any Lipschitz function interpolating the data can take. The midpoint of this range is the optimal choice in terms of minimising possible error of approximation, the norm $||f-g||_{C(X)}$.

\subsection{Monotonicity constraints}

We consider monotonicity of $f$ with respect to each variable $i$ individually, so a function can be monotone increasing in some variables, neither in others, and monotone decreasing in the rest. W.l.o.g. consider monotone increasing functions. Then there are natural tighter upper and lower bounds on monotone Lipschitz functions, namely $g(x)\leq y_k$ whenever $x \prec x^k$, and $g(x) \geq y^k$ whenever $x^k \prec x$. These extra conditions modify the distances used in \eqref{upperlower}, namely
\begin{equation}
    \label{eq:distmon}
    dist(x,x^k)=||(x-x^k)_+||
\end{equation}

Again, the function $g$ is the optimal interpolant to monotone data in the norm $||f-g||_{C(X)}$. Naturally, one expects that the data set itself is consistent with the Lipschitz and monotonicity conditions. If not the data cannot be interpolated by a monotone function, but can be approximated.

\subsection{Monotonicity on parts of domain} \label{sec:mondomain}

Imposing global monotonicity on the approximation can be too restrictive in practice (consider $f(x)=x^2$ on $[-1,1]$). For this reason we treat monotonicity on intervals as follows. Partition the domain with respect to each variable $i$ into the intervals 
\begin{eqnarray*}
    A: x_i\leq a_i \\
    AB: a_i\leq x_i \leq b_i \\ 
    B: x_i\geq b_i.
\end{eqnarray*}

\begin{figure}[!htb] \centering
\includegraphics[width=0.48\textwidth,trim=4cm 18cm 6cm 0.5cm,
    clip]{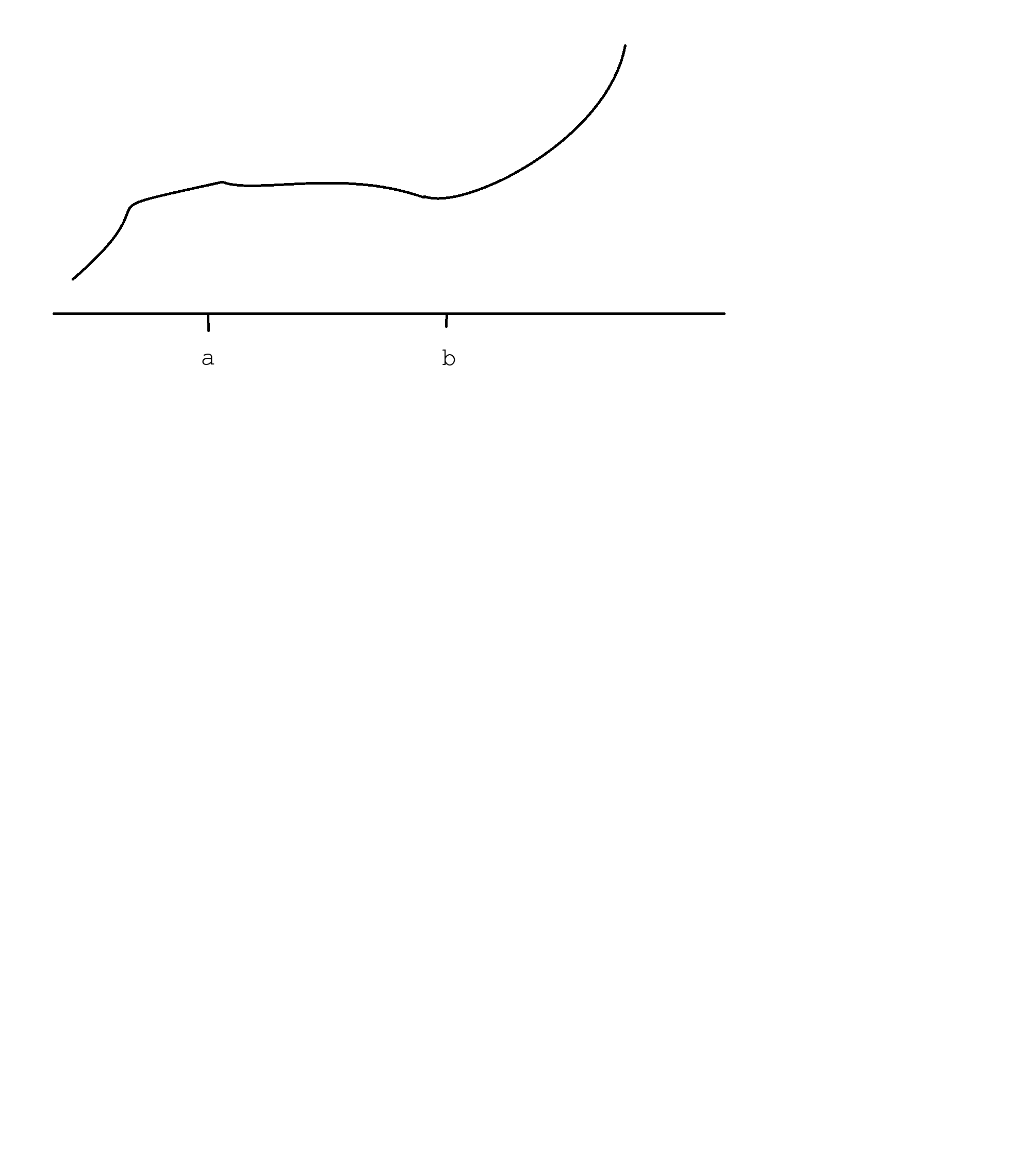}
\includegraphics[width=0.48\textwidth,trim=4cm 18cm 6cm 0.5cm,
    clip]{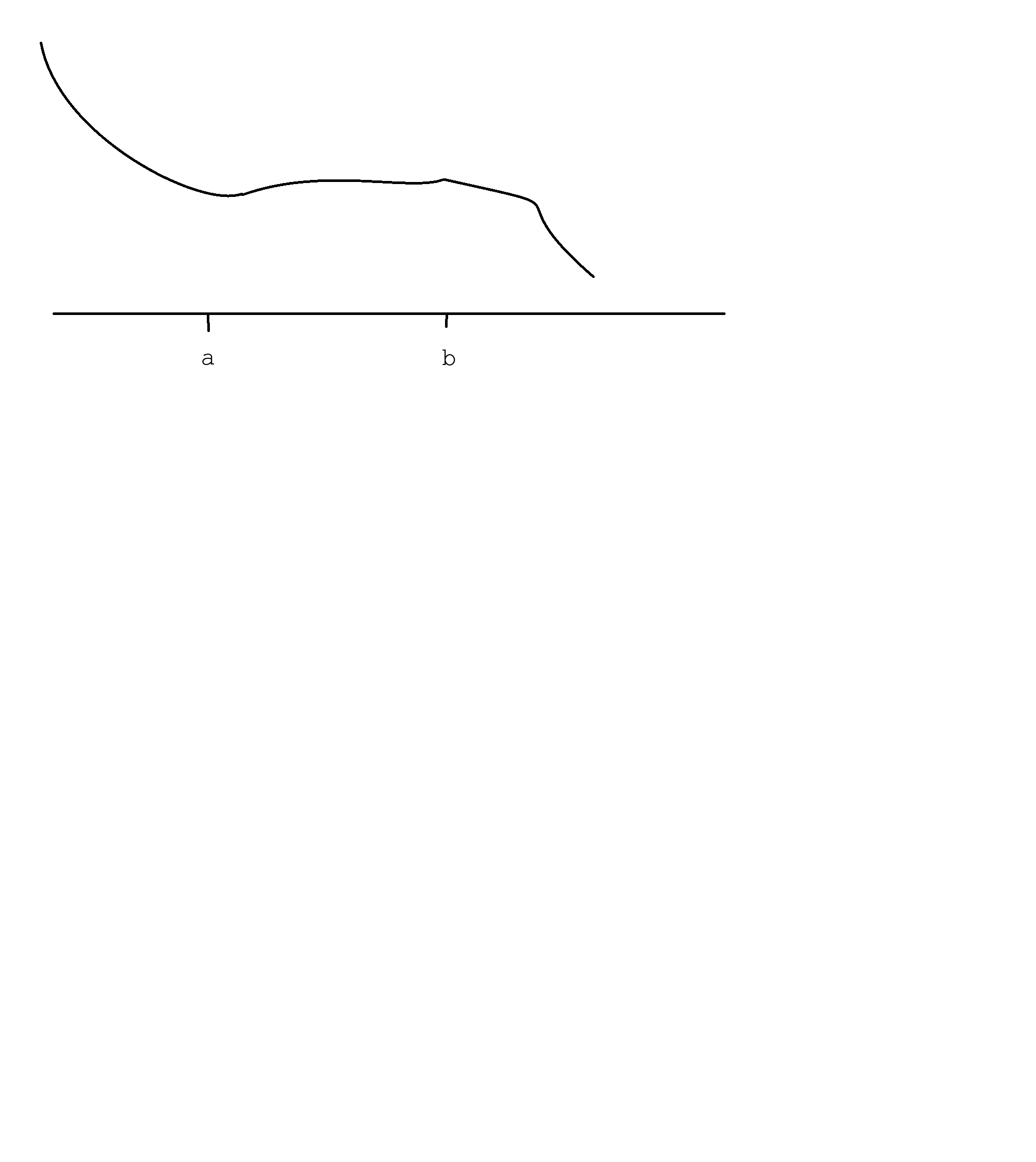}\\
\includegraphics[width=0.48\textwidth,trim=4cm 18cm 6cm 0.5cm,
    clip]{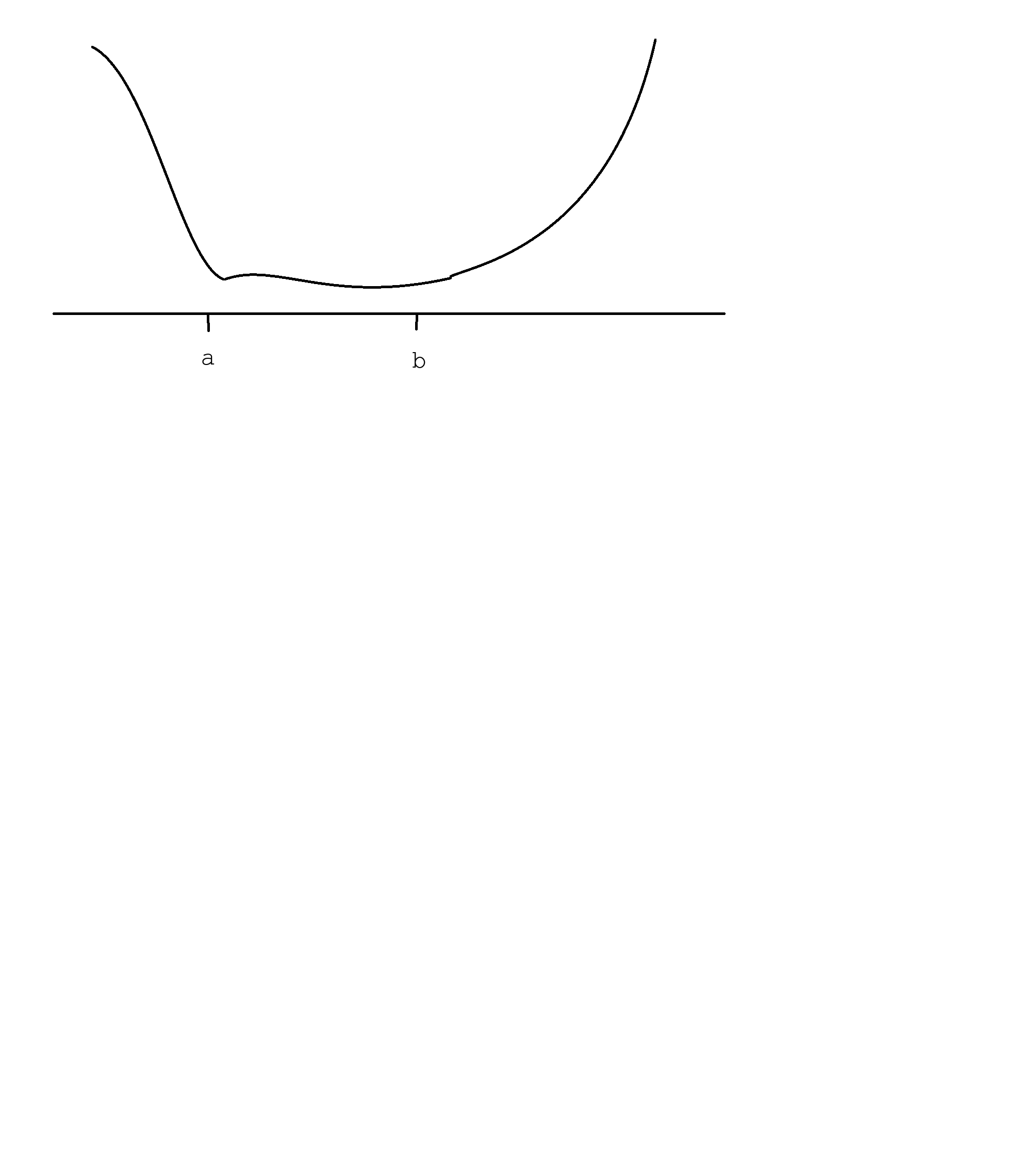}
\includegraphics[width=0.48\textwidth,trim=4cm 18cm 6cm 0.5cm,
    clip]{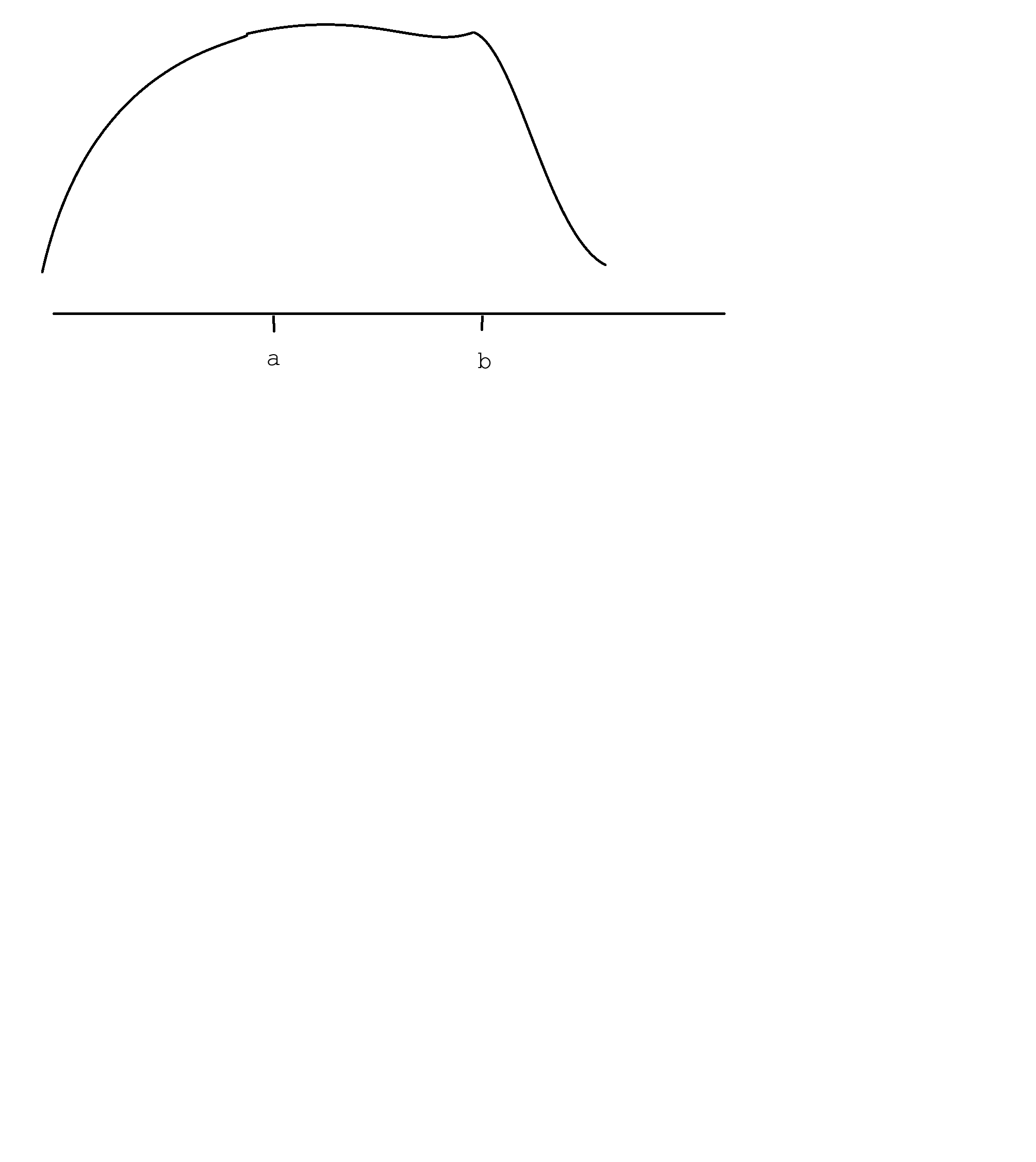}    
\caption{Four cases of restricted monotonicity on parts of the domain (cases I, II in the top row, cases III,IV in the bottom row).}\label{fig:n6}
\end{figure}

Consider four cases (see Fig. \ref{fig:n6})
$$
\mbox{Case I: } f \mbox{ increasing on }A,B, \mbox{ unrestricted on }AB.
$$
$$
\mbox{Case II: } f \mbox{ decreasing on }A,B, \mbox{ unrestricted on }AB.
$$
$$
\mbox{Case III: } f \mbox{ decreasing on }A, \mbox{ increasing on }B, \mbox{ unrestricted on }AB.
$$
$$
\mbox{Case IV: } f \mbox{ increasing on }A, \mbox{ decreasing on }B, \mbox{ unrestricted on }AB.
$$

The modification to the upper and lower bounds involve the values $a_i,b_i$ and can be expressed through components of the distance vector.

\subsection{Local Lipschitz approximation}

The Lipschitz constant $M$ used in the definition of upper and lower bounds has to be sufficiently large to ensure data can be interpolated. For some functions which have different Lipschitz properties on different parts of their domain the global value $M$ can produce rugged approximations (e.g., very large choices of $M$ produce nearly piecewise constant approximations). Instead of using one global value $M$, local Lipschitz interpolation adjusts $M$ according to the position $x$, i.e., used a function $M(x)$ in calculating the bounds. Local Lipschitz functions and the resulting geodesic distances were treated in \cite{Beliakov2007-zd}, and here we apply an approximation to $M(x)$ based on radial distances from data $x^k$, computed from the data in the form of piecewise linear spline. 

Namely, for every datum $x^k$ we calculate an approximation $M_k(||x-x^k||)$ as a non-decreasing step function. The first segment $t\in[0,d_{(1)}]$ has the value $M(t)=|y_k-y_{(1)}|/d_{(1)}$, where $d_{(1)}$ is the distance from $x^k$ to its closest neighbour, call it $t_1=d_{(1)}$. Clearly this is the smallest value consistent with the data. Next, on 
$t\in[0,t_2]$, $t_2 \geq d_{(2)}$, $M(t)=\max_j |y_k-y_j|/d_j $, and $j$ runs over the indices of the data at most at distance $t_2$ from $x^k$.  As for spline knots $t_j$ we use quantiles of the  distribution of distances from a random point $x_k$, approximated from a sample of interpoint distances. 

After this preprocessing step of building spline approximations to $M(x)$, for each datum, i.e. functions $M_k(t)$, we evaluate $g$ at query points using \eqref{upperlower} with global constant $M$ replaced with $M_k(||x-x^k||)$. The resulting approximation interpolates the data and does not result in unnecessarily large slopes in the regions where the values $y_k$ do not change much.

\subsection{Smoothing}

The smoothing method from \cite{Beliakov2007-zd} is based on minimising the sum of absolute changes to values $y_k$ that result in the interpolant consistent with any desired Lipschitz constant $M$ provided as a parameter. This problem is formulated as a linear programming problem (with $N^2$ variables and constraints), solved by the simplex method. It is efficient for up to $N=1000$. Linear programming is not particularly suitable for GPU because of its sequential  nature. 

An alternative more suitable for GPU is to transfer the following constraints for all $i \neq k$
\begin{eqnarray}
\label{eq:const1}    y_i + \varepsilon_i \leq y_k + \varepsilon_k + M d_{ik} \\
\label{eq:const2}        y_i + \varepsilon_i \geq y_k + \varepsilon_k - M d_{ik}
\end{eqnarray}
into the objective and minimise constraint violations (with penalty term $p>0$)
\begin{eqnarray} \label{eq:optim}
\nonumber    \mbox{ Min } \frac{1}{p}||\varepsilon||^2 + 
    \sum_{i\neq k}\left( (y_i-y_k - M d_{ik} -(\varepsilon_k-\varepsilon_i))_+ \right)^2\\+ \sum_{i\neq k}\left( (y_k-y_i - M d_{ik} +\varepsilon_k-\varepsilon_i)_+\right)^2.
\end{eqnarray}

Here we can use an unconstrained minimisation method like LBFGS, noticing that each objective function evaluation can be done in parallel on GPU in one step. The values $y_k-y_i$ and $d_{ik}$ are precomputed also in parallel.

\section{LipFit GPU-friendly package}

There is a lot of parallelism in all formulas presented so far, which prompts parallel Single Instruction Multiple Data (SIMD) type implementation suitable for GPUs. In this work we used Python Tensor package which provides the necessary primitives for most calculations.

We start with Lipschitz interpolation formulas \eqref{upperlower}. Clearly for every query point $x$ evaluation of distances in parallel is straightforward using torch function \texttt{torch.cdist}, and then minimum and maximum are evaluated with \texttt{torch.topk}. Moreover, multiple queries can also be executed in parallel in blocks rather than one by one, as all operations follow SIMD scenario.
Monotonicity is handled in a similar way, with a suitable modification to distance calculation, which avoids branching.

Calculation of the global Lipschitz constant from data involves calculation of all pairwise distances and then the ratios $\frac{|y_i-y_k|}{dist(x^i,x^k)}$, which is again embarrassingly parallel. There are $N^2$ operations, but no storage of the distance matrix is required. 

Next, local Lipschitz constants and interpolant evaluation are also done in parallel, with $M$ replaced with $M_k(dist(x,x^k))$. Here the primitive \texttt{torch.searchsorted} is used internally to find the required spline value, in addition to the arithmetic primitives.

The Lipschitz smoothing is emplemented through solving \eqref{eq:optim} by the internal \texttt{torch.optim.LBFGS} function, with parallel evaluation of the objective after precomputing $d_{ik}$ and $y_k-y_i$ using \texttt{torch.sum} and \texttt{torch.relu} primitives. However one needs to exercise care for large data sets $N>50,000$, as typically GPU would be unable to store several $N\times N$ matrices. For instance when $N=100,000$, a 40 GB GPU memory would be required to store $d_{ik}$.

An alternative in this case is data preprocessing, which eliminates constraints from \eqref{eq:const1},\eqref{eq:const2} which are likely to be inactive. For example, if condition
$$
y^k-\frac{M}{3}d_{ik}\leq y_i \leq y^k+\frac{M}{3} d_{ik}
$$
holds, then it is unlikely constraints \eqref{eq:const1},\eqref{eq:const2} will be violated when $\varepsilon_i, \varepsilon_k$ are added, hence such pairs $i,k$ can be removed. This is a heuristic which can save on RAM requirements, but there are cases when it does not deliver benefits (e.g. when data are generated by a linear function with some noise).

The Python package LipFit operates in the following way. First, a LipFit object (class) is created using \newline 
\texttt{LL =  LipFit(reservedata, dim, k=1, device)}, 
\newline where the device is typically "cuda" or "mps" on Mac. This class reserves space to add data dynamically on GPU, without extra copies. Then the data set is added ( and hence copied into GPU memory) with \texttt{LL.add(x,y)}, where $x \in \Re^{N \times dim}$ and $y \in \Re^N$. After that one can evaluate the interpolant at query points with \texttt{yp=LL.values(query, LC, 0, 1)} where $query \in \Re^{N_q \times dim}$ and $LC$ is the assumed Lipschitz constant (or its overestimate). The other parameters are $model=0$ (no monotonicity assumed) and $k=1$, which indicates one best upper bound and one best lower bounds are used in \eqref{upperlower}.

The new data can be added at any time and the function \texttt{values} will automatically use them for evaluation. Lipschitz constant can be found from the data using \texttt{LC=LL.lipschitz\_constant()}. The data is removed with \texttt{LL.clear()}.

When monotonicity is assumed, then one should first set monotonicity parameters \texttt{LL.setparams(mon,a,b,w)}, where for simple monotonicity, $mon\in \{-1,0,1\}^{dim}$ indicates monotonicity condition for each variable, and the rest can be \texttt{None}. When monotonicity is restricted to regions, one needs to specify $a,b\in \Re^{dim}$ and type of monotonicity $mon\in\{0,1,2,3,4\}^{dim}$ according to four cases in Section \ref{sec:mondomain}. After that evaluation is performed using \texttt{value} with $model=1$ (simple monotonicity) or $model=2$ (restricted domain monotonicity).

Calculation of local Lipschitz conditions is performed first with \texttt{LL.compute\_local\_lipschitz()} and then using \texttt{LL.values\_local(query, model, k)}. Same convention as to monotonicity types as above.

Smoothing is performed using \texttt{LL.smooth\_lipschitz( LC, model=0)}, where $LC$ is the desired Lipschitz constant, and $model$ can be 0 or 1 (with or without monotonicity). After that call, the internal data $y$ is modified with $y+\varepsilon$, and one can immediately evaluate the interpolant using \texttt{value} with the new Lipschitz constant.

When new data is added to $LL$, naturally the calculation of the Lipschitz constant, smoothing or calculation of local Lipschitz constant need to be repeated.

\paragraph{Toy example}

\begin{verbatim}
import torch
from lipfit import LipFit
import numpy as np

dim=3
ndata=1000
ntest=1000
reservedata=1000000
LL =  LipFit(reservedata,dim, 1, device)
x = torch.rand(ndata, dim)
y = torch.zeros(ndata)
xt = torch.rand(ntest, dim)

for i in range(0,ndata):
    y[i]=testfunction(x[i])
LL.add(x,y)

LC=LL.lipschitz_constant() 
mon=torch.ones(dim, dtype=torch.int32)
LL.setparams(mon,None,None,None)
yp=LL.values(xt, LC, model=0, 1)
ypmonotone=LL.values(xt, LC, model=1, 1)
\end{verbatim}

\paragraph{Extra features}

\begin{enumerate}
    \item The Lipchitz constant can be specified as a smaller value than the true Lipschitz constant, without performing Smoothing (if it is too expensive). In that case the data will be approximated as best as it can, but not in the optimal way as by \texttt{smooth\_lipschitz}.
    \item Another way of using approximation is to set $k>1$. In this case a weighting vector $w\in \Re_+^k$ should be supplied via \texttt{setparams}. The values in $w$ are assumed to add to 1, these are the weights of the "nearest" neighbour, second nearest and so on, so that both upper and lower bound are calculated using $k$ neighbours. The data set will not be interpolated. For example one can set $w=(0.8,0.1,0.1)$ and $k=3$. That $k$ has to be used in \texttt{values} and \texttt{values\_local}.
\end{enumerate}

\section{Features of instance-based approximation methods}

One can classify methods of data approximation into two broad categories: model-based and instance-based. Model-based approaches include classical linear and nonlinear regression, spline approximation, decision trees, neural networks, and similar methods which rely on a usually expensive preprocessing (training, learning) step which produces a vector of optimal parameters (coefficients) for the chosen model. Its main advantage is a compact model, inexpensive to evaluate at query points, and capturing global behaviour of the underlying function (e.g. trends, suitable for extrapolation). Saying so, models with a large number of parameters (e.g., deep learning neural networks) are typically not small, they may have large vectors of parameters comparable in size to the training data. 

On the other hand, instance-based methods have no underlying model and use all available data to make predictions at query points based on similar data. Here the method of k nearest neighbours (kNN) is a prototypical example. There is no training phase (hence lazy learning) but the evaluation becomes more expensive due to the need to process all the data. These methods can be efficiently parallelised which makes the evaluation speed less of an issue on GPUs. Another useful feature is inexpensive addition or removal of data instances, which makes them suitable for dynamic environments.

Lipschitz approximation falls into the second category, although it does not suffer from discontinuity of the approximation in kNN methods. It has several desirable features, such as preservation of the range of the data, Lipschitz continuity, optimality in the worst case scenario (adversary behaviour of the underlying function), and can also preserve monotonicity of the data (or assumed monotonicity of the underlying function). 

Model-based and instance-based approaches are complementary, both have their ranges of applicability.
Table \ref{t:2} lists some of the features of instance-based lazy learning compared to model-based prediction.

\section{Conclusions}

Lipschitz interpolation and smoothing offer  a range of desirable features of an instance-based approximation model. It provides predictable function behaviour even in adverse cases, and can incorporate additional monotonicity constraints of different types. A GPU-friendly package LipFit is available for python from \url{https://pypi.org/search/?q=lipfit}

\begin{table}
\begin{center}
\begin{tabularx}{\textwidth}{X*{4}{>{\footnotesize\arraybackslash}X}} 
Feature  & instance-based &  model-based  \\
 \hline
 Philosophy & rely on data directly (plus a few generic properties): predictions are obtained from similar observations & rely on an underlying model or function class to approximate  the behaviour\\ \\
Training & Lazy learning: no explicit training phase (optional: data smoothing and thinning) & training/fitting is essential to learn parameters  \\\\
Data dynamics & eay to add and remove data points at any time at negligible cost& expensive to update: requires re-training or fine-tuning the model \\\\

Behaviour & driven by the local structure of the data, captures complex nonlinear patterns, robust in terms of range, respects  monotonicity constraints if present in the data & learned global mapping: predictable for parametric models, can be unpredictable for black-box models, may violate known constraints \\\\
Interpretability & high (local): prediction at a query point is based on the similar observations & typically low in black-box model, higher in parametric models (linear, GAMs)\\\\
Smoothness & piecewise smooth Lipschitz-continuous, differentiable almost everywhere & usually smooth by construction, may result in large Lipschitz constant \\\\
Locality & local by design, predictions depend on similar data & global: predictions may depend on training data away from query point\\\\

Evaluation of queries & cost grows with data size, highly parallelisable (excellent for  GPU) & fast for small models, slower when model size is large (comparable to  training size as in many deep learning models)\\
\hline
\end{tabularx}
\caption{Features and drawbacks of instance-based and model-based approximation} \label{t:2}
\end{center}
\end{table}
\begin{table}
\begin{center}
\begin{tabularx}{\textwidth}{X*{4}{>{\footnotesize\arraybackslash}X}} 
Feature  & instance-based &  model-based  \\
 \hline
Hyperparameter sensitivity& sensitive to the choice of $k$, distance metric and kernel & sensitive to model architecture, regularisation and learning rate\\\\
Robustness to noise& robust with appropriate weighting, local influence limits damage& varies: some models are sensitive unless regularised\\\\
Scalability & suffers from the curse of dimensionality& also affected but low rank/sparse model scale better\\\\
Memory requirements& high: stores all data& low for models with few parameters, high for many deep learning models \\\\
Extrapolation& typically poor outside the convex hull of data& can extrapolate (quality depends heavily on model bias and assumptions)
\\
\hline
\end{tabularx}
\caption{Table continued}
\end{center}
\end{table}

\newpage
 \bibliographystyle{abbrvurl}
\bibliography{sample}

\end{document}